\documentclass[12pt]{article}
\usepackage{amsmath,amssymb,amsthm, amscd}

\begin{document}
\begin{center} 
{\bf Poisson deformations and birational geometry}\footnote{This is a revised version 
of the article entitled ``Poisson deformations and Mori dream spaces".}  
\end{center} 
\vspace{0.4cm}

\begin{center}
Yoshinori Namikawa
\end{center} 
\vspace{0.4cm}

\begin{center}
{\bf \S 1. Introduction}  
\end{center}

A complex normal variety $X$ is a {\em symplectic variety} if there is a holomorphic 
symplectic 2-form $\omega$ on the regular locus $X_{reg}$ of $X$ and $\omega$ 
extends to a holomorphic 2-form on a resolution $\tilde{X}$ of $X$. 
Notice that, if $\omega$ extends to a holomorphic 2-form on a particular resolution 
of $X$, then it also extends to a holomorphic 2-form on an arbitrary resolution of 
$X$.      

In this article $X$ is an affine symplectic variety with a coordinate ring $R$, where $R$ is 
positively graded as $R = \oplus_{i \geq 0} R_i$ with $R_0 = \mathbf{C}$. Moreover 
we assume that $\omega$ is homogeneous with respect to the natural $\mathbf{C}^*$-action 
on $X$. The weight $l$ of $\omega$ is automatically positive because $X$ has only 
canonical singularities ([Na 1], Lemma 2.2). Such varieties are constructed in various ways 
such as nilpotent orbit closures of a complex semisimple Lie algebra, Slodowy slices to 
nilpotent orbits, symplectic quotient singularities, Nakajima quiver varieties and so on. 

In the remainder we assume that $X$ has a (projective) crepant resolution $\pi: Y 
\to X$\footnote{More generally we may take a {\bf Q}-factorial terminalization as $Y$, 
which always exists by a main result of [BCHM].}  
Then the symplectic 2-form 
$\omega$ extends to a symplectic 2-form $\omega_Y$ on $Y$. Moreover the 
$\mathbf{C}^*$-action uniquely extends to a $\mathbf{C}^*$-action on $Y$ ([Na 2], Proposition A.7).  Notice that the symplectic structures $\omega$ and 
$\omega_Y$ determine Poisson structures respectively on $X$ and $Y$. 
A Poisson deformation of $X$ (resp. $Y$) is just a deformation of the 
pair of $X$ (resp. $Y$) and its Poisson structure. These have been extensively studied 
in [Ka],  [K-V], [G-K], [Na 3] and [Na 4].  

Here let us consider the second cohomology space $H^2(Y, \mathbf{C})$. 
The first interpretation of this space is the Picard group 
$\mathrm{Pic}(Y) \otimes_{\mathbf Z}\mathbf{C}$ tensorised with $\mathbf{C}$. 
One can consider various cones inside $\mathrm{Pic}(Y) \otimes_{\mathbf Z}\mathbf{R}$. 
In particular, the $\pi$-movable cone $\mathrm{Mov}(\pi)$ is decomposed into 
the ample cones $\mathrm{Amp}(\pi')$ for different crepant resolutions $\pi': Y' \to X$ 
of $X$. The codimension one faces of each ample cone generate hyperplanes of  
$\mathrm{Pic}(Y) \otimes_{\mathbf Z} \mathbf{R}$ and these hyperplanes divide 
$\mathrm{Mov}(\pi)$ into various chambers. 

The second interpretation of $H^2(Y, \mathbf{C})$ is the {\em base space} of the universal Poisson deformation of $Y$.   By using the $\mathbf{C}^*$-action on $Y$ one 
may construct the universal Poisson deformation $\mathcal{Y} \to H^2(Y, \mathbf{C})$. 
The central fibre $\mathcal{Y}_0$ over $0 \in H^2(Y, \mathbf{C})$ is isomorphic to $Y$; 
hence $\mathcal{Y}_0$ usually contains proper subvarieties as the exceptional locus of 
$\pi$. But it turns out that a general fibre $\mathcal{Y}_t$ is an affine variety. 
We define the subset $\mathcal{D} \subset H^2(Y, \mathbf{C})$ as the locus where 
$\mathcal{Y}_t$ is not affine. 

The purpose of this article is to compare these two structures on $H^2(Y, \mathbf{C})$. 
We explain the main result in a more precise form.   
The universal Poisson deformation of $X$ is constructed as a family of Poisson varieties 
over a quotient space $H^2(Y, \mathbf{C})/W$ by a finite group $W$. The universal Poisson 
deformations of $X$ and $Y$ fit into the commutative diagram

\begin{equation} 
\begin{CD} 
\mathcal{Y} @>>> \mathcal{X} \\ 
@VVV @VVV \\ 
H^2(Y, \mathbf{C}) @>{q}>> H^2(Y, \mathbf{C})/W      
\end{CD} 
\end{equation}
 
Here the central fibres $\mathcal{Y}_0$ and $\mathcal{X}_{q(0)}$ are respectively 
$Y$ and $X$ and the induced map $\mathcal{Y}_0 \to \mathcal{X}_{q(0)}$ coincides 
with $\pi$.  

Put $\mathcal{X}' := \mathcal{X} \times_{H^2(Y, \mathbf{C})/W} H^2(Y, \mathbf{C})$ and 
consider the induced commutative diagram 

\begin{equation} 
\begin{CD} 
\mathcal{Y} @>{\Pi}>> \mathcal{X}' \\ 
@VVV @VVV \\ 
H^2(Y, \mathbf{C}) @>{id}>> H^2(Y, \mathbf{C})     
\end{CD} 
\end{equation}

The maps $\Pi_t: \mathcal{Y}_t \to \mathcal{X}'_t$ are birational for 
all $t \in H^2(Y, \mathbf{C})$ and  
they are isomorphisms for general $t \in H^2(Y, \mathbf{C})$. 
Let $\mathcal{D} \subset H^2(Y, \mathbf{C})$ be the locus where 
$\Pi_t$ is not an isomorphism. Notice that $\mathcal{D}$ coincides with 
the locus where $\mathcal{X}'_t$ is singular. It can be also characterized as the 
locus where $\mathcal{Y}_t$ is not an affine variety.  
Then our main theorem asserts that
\vspace{0.15cm}

{\bf Main Theorem}:  
\vspace{0.15cm} 

(i) there are finitely many linear spaces $\{H_i\}_{i \in I}$ 
of $H^2(Y, \mathbf{Q})$ with $\mathrm{Codim}(H_i) = 1$ such that $\mathcal{D} 
= \bigcup_{i \in I} (H_i)_{\mathbf C}$, 

(ii) there are only finitely many crepant projective resolutions $\{\pi_k : Y_k \to X\}_{k \in K}$ of $X$, and 

(iii) the set of open chambers determined by the real hyperplanes 
$\{(H_i)_{\mathbf R}\}_{i \in I}$ coincides with the set $\{w(\mathrm{Amp}(\pi_k))\}$ 
where $w \in W$ and $\pi_k$ are crepant projective resolutions of $X$.     
\vspace{0.2cm}

Notice that a crepant resolution $\pi: Y \to X$ is a relative Mori dream space in the sense of 
[H-K] when $X$ is more generally a rational Gorenstein singularities (cf. Corollary 1.3.2 of [BCHM]).  But the chamber structure is determined by the set of {\em full} hyperplanes when $X$ is a symplectic singularity. This property does not hold for a general rational Gorenstein  singularity (Example 13).

The simplest nontrivial example of the main theorem would be the following. 
\vspace{0.2cm}

{\bf Example}. Let $X :=\{(x,y,z) \in \mathbf{C}^3; x^2 + y^2 + z^n = 0\}$ be an 
$A_{n-1}$-surface singularity and let $\pi: Y \to X$ be the minimal resolution. Let 
$E_1, ..., E_{n-1}$ be the irreducible components of $\mathrm{Exc}(\pi)$. 
The 2-nd cohomology $H^2(Y, \mathbf{R})$ has dimension $n-1$ and is generated
by $[E_i]$'s. We introduce a negative definite symmetric form $(\;, \;)$ on 
$H^2(Y, \mathbf{R})$ by the usual intersection pairings of $E_i$'s. 
The subset $\Phi := \{E := \Sigma a_iE_i; E^2 = -2, a_i \in \mathbf{Z}\}$ of $H^2(Y, \mathbf{R})$ determines a root system of $A_{n-1}$-type. The Weyl group of $\Phi$ 
is isomorphic to the symmetric group $S_n$. 
Thus $S_n$ acts on $H^2(Y, \mathbf{R})$ as the Weyl group. 
Moreover $H^2(Y, \mathbf{R})$ is divided into $n !$ Weyl chambers. 
In particular, $\mathrm{Amp}(\pi)$ coincides with one of them:  
$$\mathrm{Amp}(\pi) = \{x \in H^2(Y, \mathbf{R}); (x, [E_j]) > 0 \; \mathrm{for} \; \mathrm{all}\; j\}.$$ 
 
Note that $X$ has a $\mathbf{C}^*$-action $(x,y,z) \to (t^nx, t^ny, t^2z)$, $t \in 
\mathbf{C}^*$ and the symplectic 2-form 
$$\omega := \mathrm{Res}(dx \wedge dy \wedge dz/x^2 + y^2 + z^n)$$ is homogeneous 
of weight 2 with respect to this $\mathbf{C}^*$-action.   
The universal Poisson deformation of $X$ is given by 
$$ \mathcal{X} := \{(x, y, z, u_1, ... u_{n-1}) \in \mathbf{C}^{n+2}; x^2 + y^2 + z^n + 
u_1z^{n-2} + ... + u_{n-2}z + u_{n-1} = 0 \}$$ 
where $(u_1, ..., u_{n-1}) \in \mathbf{C}^{n-1}$ is the base space. Let us consider an $n-1$   dimensional vector space $V := \{(s_1, ..., s_n) \in \mathbf{C}^n; \Sigma s_i = 0\}$ and a finite 
Galios covering $V \to \mathbf{C}^{n-1}$ defined by $(s_1, ..., s_n) \to 
(\sigma_2(s_1, ..., s_n), ..., \sigma_n(s_1, ..., s_n))$.  Here $\sigma_i$ are defined by 
$$(z-s_1)(z-s_2)\cdot \cdot \cdot (z-s_n) = z^n + \sigma_1(s_1, ..., s_n)z^{n-1} + ... 
+ \sigma_n(s_1, ..., s_n).$$ 
The Galois group of this map is the symmetric group $S_n$ acting on $V$ by 
the permutation of $s_i$'s. Thus $\mathbf{C}^{n-1} = V/S_n$. 
We put $\mathcal{X}' := \mathcal{X} \times_{V/S_n} V$. Then $\mathcal{X}'$ can be 
written as 
$$ \mathcal{X}' := \{(x, y, z, s_1, ..., s_n) \in \mathbf{C}^{n+3}; x^2 + y^2 + 
(z - s_1)(z - s_2)\cdot\cdot\cdot (z-s_n)  = 0, \; \Sigma s_i = 0\}.$$ 
 
The map $\mathcal{X}' \to V$ can be viewed as an $n-1$ dimensional family of surfaces. 
Define the discriminant locus $\mathcal{D}$ by  
$$\mathcal{D} := \{s = (s_1, ..., s_n) \in V; \mathcal{X}'_s \; \mathrm{is} \; \mathrm{singular} 
\}.$$  Then one has  
$$\mathcal{D} = \bigcup_{1 \le i < j \le n} L_{ij},$$ where $$L_{ij} = \{(s_1, ..., s_n) \in V; s_i = s_j\}.$$
The family can be resolved simultaneously. The simultaneous resolution is not unique. 
Let us take one of them, say $\mathcal{Y}$.  
Then  we have a commutative diagram    

\begin{equation} 
\begin{CD} 
\mathcal{Y} @>{\Pi}>> \mathcal{X}' \\ 
@VVV @VVV \\ 
V @>{id}>> V     
\end{CD} 
\end{equation}

Then $\mathcal{Y} \to V$ 
is the universal Poisson deformation of $Y$.    
   
For a general point $t_{ij} \in L_{ij}$, the 
induced map $\Pi_{t_{ij}}: \mathcal{Y}_{t_{ij}} \to \mathcal{X}'_{t_{ij}}$ contracts one smooth rational curve $C_{ij}$ to a point. By using the $\mathbf{C}^*$-action on $\mathcal{Y}$ 
we put $E_{ij}: =  \lim_{\lambda \to 0} \lambda(C_{ij})$  
We may choose a simultaneous resolution 
$\Pi: \mathcal{Y} \to \mathcal{X}'$ so that  
$E_{12} = E_1$, $E_{23} = E_2$, ..., $E_{n-1,n} = E_{n-1}$.     

We introduce a positive definite symmetric form on $$V(\mathbf{R}) 
:= \{(s_1, ..., s_n) \in \mathbf{R}^n; \Sigma s_i = 0\}$$ as the restriction 
of the standard metric on $\mathbf{R}^n$. 
Define an $\mathbf{R}$-linear isomorphism $\alpha : H^2(Y, \mathbf{R}) 
\to V(\mathbf{R})$ by $\alpha ([E_1]) = (1,-1,0, ..., 0)$, $\alpha ([E_2]) = 
(0,1,-1,0,...0)$, ..., $\alpha ([E_{n-1}]) = (0,0, ..., 1, -1)$. The map $\alpha$ 
preserves the symmetric forms up to a reversal of sign. 
Each Weyl chamber of $H^2(Y, \mathbf{R})$ is mapped to a chamber of 
$V(\mathbf{R})$ determined by $L_{ij} \cap V(\mathbf{R})$'s. 
In particular, $\mathrm{Amp}(\pi)$ is mapped to a chamber surrounded  by $n-1$ rays 
$$\mathbf{R}_{\geq 0} (-n+i, ..., -n+i, i, ..., i)  \;\; (1 \le i \le n-1)$$ 
Here $-n+i$ occurs $i$ times and  $i$ occurs $n-i$ times.  

The map $\alpha$ induces a $\mathbf{C}$-linear isomorphism 
$\alpha_{\mathbf C} : H^2(Y, \mathbf{C}) \to V$, whose inverse coincides with the 
period map for $\mathcal{Y} \to V$ determined by a suitable multiple of $\omega$. 
(END of EXAMPLE)   
\vspace{0.3cm}

The example above can be generalized to the Slodowy slice to an arbitrary nilpotent 
orbit (see Example 11). 

One can find in  [Na 6] explicit descriptions of the chamber structures and the flops connecting them in the case where $X$ is a nilpotent orbit closure of a complex simple Lie algebra. 
In particular, we can determine the explicit equations of the hyperplanes $H_i$. In Example 12 
we illustrate this by using one typical example.      

The birational geometry for crepant resolutions of 
4-dimensional symplectic singularities are also treated in [A-W].   
     
\vspace{0.2cm}

\begin{center}
{\bf \S 2. Proof of Main theorem}  
\end{center}

Let $\pi: Y \to X$ be the same as in Introduction. 
Let $\mathrm{Amp}(\pi) \subset H^2(Y, \mathbf{R})$ be the open convex cone 
generated by $\pi$-ample line bundles on $Y$. Its closure $\overline{\mathrm{Amp}(\pi)}$ 
coincides with the cone generated by $\pi$-nef line bundles. Thus we often call 
$\overline{\mathrm{Amp}(\pi)}$ the {\em nef cone} of $\pi$.  
A line bundle $L$ on $Y$ is called $\pi$-movable if the support of 
$\mathrm{Coker}[\pi^*\pi_* L \to L]$ has codimension at least two. 
The closed cone in $H^2(Y, \mathbf{R})$ generated by the classes of $\pi$-movable line bundles is denoted by $\overline{\mathrm{Mov}(\pi)}$ and its interior is called the 
{\em movable cone} for $\pi$.  In the remainder we will write $\mathrm{Mov}(\pi)$ 
for the movable cone.  By the definition we have $\mathrm{Amp}(\pi) \subset \mathrm{Mov}(\pi)$.

The following is a 
fundamental fact on birational geometry proved almost thirty years ago. 
\vspace{0.2cm}

{\bf Lemma 1}. {\em The nef cone $\overline{\mathrm{Amp}(\pi)}$ is a  
rational polyhedral cone.} 
\vspace{0.2cm}

{\em Proof}. Take an effective divisor $D$ of $Y$ in such a way that $-D$ is 
$\pi$-ample. If we take a rational number $\epsilon > 0$ very small and put 
$\Delta := \epsilon D$. Then $(Y, \Delta)$ has klt (Kawamata log terminal) singularities. 
Since $K_Y$ is trivial, $-(K_Y + \Delta)$ is $\pi$-ample. The cone of effective 
1-cycle cone $\overline{NE(\pi)}$ is then a rational polyhedral cone (cf. 
[KMM, Theorem 4-2-1, Proposition 4-2-4]). 
Its dual cone $\overline{\mathrm{Amp}(\pi)}$ is hence a rational polyhedral 
cone. Q.E.D.   
\vspace{0.2cm}

The symplectic 2-form $\omega$ on $X_{reg}$ identifies the sheaf $\Theta_{X_{reg}}$ 
of holomorphic vector fields with the sheaf $\Omega^1_{X_{reg}}$ of holomorphic 
1-forms. Moreover it induces an isomorphism $\wedge^2 \Theta_{X_{reg}} \cong 
\Omega^2_{X_{reg}}$. By this isomorphism $\omega$ is identified with a 2-vector 
$\theta$ on $X_{reg}$, which we call the {\em Poisson 2-vector}. 
The Poisson 2-vector defines a bracket on the structure sheaf $\mathcal{O}_{X_{reg}}$ 
by $\{f, g\} := \theta (df \wedge dg)$ with $f \in \mathcal{O}_{X_{reg}}$ and $g \in 
\mathcal{O}_{X_{reg}}$. By the definition this bracket is a $\mathbf{C}$-bilinear form 
and is a biderivation with respect to each factor.  One can prove that this 
bracket satisfies the Jacobi identity by using the fact that $\omega$ is d-closed. 
Thus $X_{reg}$ admits a {\em Poisson structure}. By the normality of $X$, this Poisson structure uniquely extends to a Poisson structure on $X$. 

A $T$-scheme $\mathcal{X} \to T$ is called a Poisson $T$-scheme if there is an $\mathcal{O}_T$-bilinear Poisson bracket $\{\;, \;\}_{\mathcal{X}} : \mathcal{O}_{\mathcal X} \times \mathcal{O}_{\mathcal X} \to \mathcal{O}_{\mathcal X}$. Let $T$ be a scheme 
over $\mathbf{C}$ and let $0 \in T$ be a closed point. 

{\em A Poisson deformation} of the 
Poisson variety $X$ is a Poisson $T$-scheme $f: \mathcal{X} \to T$ together 
with an isomorphism $\phi: \mathcal{X}_0 \cong X$  which satisfies the following conditions  

(i) $f$ is a flat surjective morphism, and 

(ii) $\{\;, \;\}_{\mathcal X}$ restricts to the original Poisson structure $\{\;, \;\}$ on 
$X$ via the identification $\phi$. 
\vspace{0.2cm} 

Two Poisson deformations $(\mathcal{X}/T, \phi)$ and $(\mathcal{X}' /T, \phi')$ with 
the same base are {\em equivalent} if there is a $T$-isomorphism $\mathcal{X} \cong 
\mathcal{X}'$ of Poisson schemes such that it induces the identity map on the 
central fibre.  

For a local Artinian $\mathbf{C}$-algebra $A$ with residue field $\mathbf{C}$, 
denote by $\mathrm{PD}_X(A)$ the set of equivalence classes of Poisson deformations of 
$X$ over $\mathrm{Spec} (A)$. Then it defines a functor   
$$\mathrm{PD}_X:  (Art)_{\mathbf C} \to (Set)$$ from the 
category of local Artinian $\mathbf{C}$-algebra with residue field $\mathbf{C}$ to 
the category of sets.   
  
Assume that a crepant projective resolution  $\pi: Y \to X$ is given. 
Then the symplectic 2-form 
$\omega$ extends to a symplectic 2-form $\omega_Y$ on $Y$.  The symplectic 2-form $\omega_Y$ determines a 
Poisson structure $\{\;, \;\}_Y$ on $Y$. Notice that $\pi: (Y, \{\;, \;\}_Y) \to (X, \{\;, \;\})$ 
is a morphism of Poisson varieties, i.e. $\pi^*\{f,g\} = \{\pi^*f, \pi^*g\}_Y$ for $f, g \in 
\mathcal{O}_X$. 

Assume that $A \in (Art)_{\mathbf C}$ and $T := \mathrm{Spec} A$.  
If $(\mathcal{Y}/T, \phi)$ is a Poisson deformation of $Y$, then $\mathrm{Spec} 
\Gamma (Y, \mathcal{O}_{\mathcal Y}) \to T$ is a flat deformation of $X = \mathrm{Spec}
\Gamma(Y, \mathcal{O}_Y)$ because $H^1(Y, \mathcal{O}_Y) = 0$ (cf. [Wa]). 
Moreover, the Poisson structure on $\mathcal{Y}$ determines a Poisson bracket on 
$\Gamma (Y, \mathcal{O}_{\mathcal Y})$. Thus we have a morphism of functors
$$ \pi_* : PD_Y \to PD_X.$$  
We can apply Schlessinger's theory [Sch] to these functors to get the prorepresentable 
hulls $R_X$ and $R_Y$. By definition $R_X$ and $R_Y$ are both complete local 
$\mathbf{C}$-algebras with residue field $\mathbf{C}$. Let $m_X$ (resp. $m_Y$) be 
the maximal ideal of $R_X$ (resp. $R_Y$) and put $R_{X,n} := R_X/m_X^{n+1}$ 
and $R_{Y, n} := R_Y/m_Y^{n+1}$. There are formal semiuniversal Poisson deformations 
$\{X_n \to \mathrm{Spec} R_{X,n}\}_{n \geq 0}$ and 
$\{Y_n \to \mathrm{Spec} R_{Y,n}\}_{n \geq 0}$. 
Moreover, for each $n \geq 0$, there is a commutative diagram 

\begin{equation} 
\begin{CD} 
Y_n @>>> X_n \\ 
@VVV @VVV \\ 
\mathrm{Spec} R_{Y,n} @>>> \mathrm{Spec} R_{X,n}     
\end{CD} 
\end{equation}

By Corollary 15 of [Na 2] the Poisson deformation functor $PD_Y$ is unobstructed 
and its tangent space $PD_Y (\mathbf{C}[\epsilon])$ is naturally isomorphic to 
$H^2(Y, \mathbf{C})$ (cf. [Na 2], Proposition 8, Corollary 10).  
By Theorem 5.1 of [Na 3], $PD_X$ is also unobstructed and $R_Y$ is a finite $R_X$-module.   
Moreover, by Corollary 2.5 of [Na 3], both functors $PD_X$ and $PD_Y$ are prorepresentable, 
in other words, $\{X_n \to \mathrm{Spec} R_{X,n}\}_{n \geq 0}$ and $\{Y_n \to \mathrm{Spec} R_{Y,n}\}_{n \geq 0}$ are respectively formal universal Poisson deformations of $X$ and $Y$.         

The main result of [Na 4] asserts that $\pi_*: PD_Y \to PD_X$ is a finite {\em Galois} 
covering with a finite Galois group $W$. 
By the $\mathbf{C}^*$-actions we can algebraize the commutative diagram above  
([Na 3, Section 5])      

\begin{equation} 
\begin{CD} 
\mathcal{Y} @>>> \mathcal{X} \\ 
@VVV @VVV \\ 
H^2(Y, \mathbf{C}) @>{q}>> H^2(Y, \mathbf{C})/W      
\end{CD} 
\end{equation}

Here $W$ acts linearly on $H^2(Y, \mathbf{Q})$ 
and the diagram is 
$\mathbf{C}^*$-equivariant.
Each linear function on $H^2(Y, \mathbf{C})$ has weight $l = wt(\omega)$.       
For each $t \in H^2(Y, \mathbf{C})$, the induced map $\mathcal{Y}_t \to \mathcal{X}_{q(t)}$ 
is a birational morphism and it is an isomorphism for general $t$. Notice that $\mathcal{Y}_0 
\to \mathcal{X}_{q(0)}$ coincides with $\pi: Y \to X$.  The map $\mathcal{Y} \to 
H^2(Y, \mathbf{C})$ is a $C^{\infty}$-fibre bundle and every fibre $\mathcal{Y}_t$ is 
diffeomorphic to $Y$. In particular, there is a natural identication of $H^2(\mathcal{Y}_t, \mathbf{Z})$ with $H^2(Y, \mathbf{Z})$. 

Put $\mathcal{X}' := \mathcal{X} \times_{H^2(Y, \mathbf{C})/W} H^2(Y, \mathbf{C})$ and 
consider the induced commutative diagram 

\begin{equation} 
\begin{CD} 
\mathcal{Y} @>{\Pi}>> \mathcal{X}' \\ 
@V{f}VV @V{g}VV \\ 
H^2(Y, \mathbf{C}) @>{id}>> H^2(Y, \mathbf{C})     
\end{CD} 
\end{equation}

Here $\Pi : \mathcal{Y} \to \mathcal{X}'$ is a crepant projective resolution. 
Let $\Pi': \mathcal{Y}' \to \mathcal{X}'$ be another crepant projective resolution. 
Since the birational map $\mathcal{Y}' --\to \mathcal{Y}$ is an isomorphism in 
codimension one, we have a natural identification $H^2(\mathcal{Y}', \mathbf{R}) 
\cong H^2(\mathcal{Y}, \mathbf{R})$. By this identification the nef cone 
$\overline{\mathrm{Amp}(\Pi')}$ can be seen as a cone in $H^2(\mathcal{Y}, \mathbf{R})$. 
Furthermore, as the restriction map $H^2(\mathcal{Y}, \mathbf{R}) \to H^2(Y, \mathbf{R})$ 
is an isomorphism, we can regard $\overline{\mathrm{Amp}(\Pi')}$ as a cone in 
$H^2(Y, \mathbf{R})$.  By the same reasoning as in Lemma 1 $\overline{\mathrm{Amp}(\Pi')}$ is a rational polyhedral cone. 
\vspace{0.2cm}  
 
Take a nonzero element $\alpha \in H^2(Y, \mathbf{C})$ and consider the complex 
line $\mathbf{A}^1_{\alpha} := \{t\alpha\}_{t \in \mathbf{C}}$ inside $H^2(Y, \mathbf{C})$. 
We put 
$$\mathcal{Y}^{\alpha} :=  \mathcal{Y} \times_{H^2(Y, \mathbf{C})} \mathbf{A}^1_{\alpha},$$ 
and 
$$\mathcal{X}^{\alpha} := \mathcal{X}' \times_{H^2(Y, \mathbf{C})} \mathbf{A}^1_{\alpha}.$$ 
We then get a commutative diagram 

\begin{equation} 
\begin{CD} 
\mathcal{Y}^{\alpha} @>{\Pi^{\alpha}}>> \mathcal{X^{\alpha}} \\ 
@VVV @VVV \\ 
\mathbf{A}^1_{\alpha} @>{id}>> \mathbf{A}^1_{\alpha}    
\end{CD} 
\end{equation}

According to Kaledin [Ka] we call $\mathcal{Y}^{\alpha} \to \mathbf{A}^1_{\alpha}$ the 
{\em twistor deformation} of $Y$ determined by $\alpha$. There is a relative symplectic 
2-form $\omega_{\mathcal Y}$ on $\mathcal{Y}$ extending the symplectic 2-form $\omega_Y$ 
on $Y$. Define $\omega_{{\mathcal Y}^{\alpha}}$ to be the restriction of $\omega_{\mathcal Y}$ 
to $\mathcal{Y}^{\alpha}$. A remarkable property of the twistor deformation is that, for 
each $t \in \mathbf{A}^1_{\alpha}$, the 2-nd cohomology class  
$[\omega_{\mathcal{Y}^{\alpha}_t}] \in H^2(\mathcal{Y}^{\alpha}_t, \mathcal{C})$ 
coincides with $t\cdot \alpha \in H^2(Y, \mathbf{C})$ under the natural identification 
$H^2(\mathcal{Y}^{\alpha}_t, \mathcal{C}) \cong H^2(Y, \mathbf{C})$ (cf. [G-K]).  
\vspace{0.2cm}

{\bf Lemma 2}. {\em If $(\alpha, C) \ne 0$ for any proper curve $C \subset Y$ 
such that $\pi (C)$ is a point, then $\Pi^{\alpha}_t$ is an isomorphism for 
all $t \in \mathbf{A}^1_{\alpha} - \{0\}$.} 
\vspace{0.2cm}

{\em Proof}. Assume that a proper curve $C_t \subset \mathcal{Y}^{\alpha}_t$ 
is contracted to a point by $\Pi^{\alpha}_t$ for some $t \in \mathbf{A}^1_{\alpha} - \{0\}$. 
Since $\omega_{\mathcal{Y}^{\alpha}_t}$ is a holomorphic 2-form, we must have 
$([\omega_{\mathcal{Y}^{\alpha}_t}], C_t) = 0$. By the $\mathbf{C}^*$-action, we have a 
proper (but, not necessarily irreducible or reduced) curve $C \subset Y$ as a limit of 
$C_t$. Then we have $(t \alpha, C) = 0$; hence $(\alpha, C) = 0$.     
Q.E.D. 
\vspace{0.2cm}

It is not true that if $(\alpha, C) = 0$ for some proper curve $C \subset Y$, then 
the birational map $\Pi^{\alpha}_t$ contracts some curves to points for each 
$t \in \mathbf{A}^1_{\alpha}$. But we have the following: 
\vspace{0.2cm}

{\bf Lemma 3}. {\em Assume that $\pi$ factorizes 
as $Y \stackrel{\mu}\to Z \to X$ with $Z$ normal, $\mu$ birational, and $\mathrm{Exc}(\mu) \ne \emptyset$.  
If  $\alpha \in \mu^*(H^2(Z, \mathbf{C}))$, then 
$\mathrm{Exc}(\Pi^{\alpha}_t) \ne \emptyset$ for all $t \in \mathbf{A}^1_{\alpha}$.}  
\vspace{0.2cm}

{\em Proof}. Since $Z$ has only rational singularities, we have an exact sequence 
$$ 0 \to H^2(Z, \mathbf{C}) \to H^2(Y, \mathbf{C}) \to H^0(Z, R^2\mu_*\mathbf{C}).$$
Recall that a 1-st order Poisson deformation $\mathcal{Y}_1 \to \mathrm{Spec} \mathbf{C}[\epsilon]$ of $Y$ corresponds to an element $\alpha$ of $H^2(Y, \mathbf{C})$. For a point $p \in Z$ take a sufficiently small analytic neighborhood $U$ of $p$. Then $\mathcal{Y}_1$ induces a 1-st order Poisson deformation of $\mu^{-1}(U)$. 
If $\alpha \in H^2(Z, \mathbf{C})$, then the induced Poisson deformation of $\mu^{-1}(U)$ 
is trivial by the exact sequence above.  

Now let us define a subfunctor $\mathrm{PD}^\mu_Y$ of $\mathrm{PD}_Y$ as 
follows. 
\vspace{0.2cm}

For $A \in (Art)_{\mathbf C}$, define $\mathrm{PD}^\mu_Y(A)$ to be 
the set of equivalence classes of the Poisson deformations 
$f_T: \mathcal{Y}_T \to T$ of $Y$ ($T := \mathrm{Spec}(A)$) with the following property:
 
(*)  For each $p \in Z$, there is a sufficiently small (analytic) neighborhood $U$ of $p$ 
such that $f_T$ induce trivial Poisson deformations of $\mu^{-1}(U)$. 
\vspace{0.2cm}
 
The functor $\mathrm{PD}^\mu_Y$ has a prorepresentable hull $R^{\mu}_Y$ in the sense of [Sch] and its tangent space $PD^{\mu}_Y(\mathbf{C}[\epsilon])$ is isomorphic to $H^2(Z, \mathbf{C})$. 
As the tangent space has topological nature one can apply the $T^1$-lifting method 
to prove that $PD^{\mu}_Y$ is unobstructed.

Recall that the twistor deformation $\mathcal{Y}^{\alpha}$ is constructed inductively as 
the objects of $\mathrm{PD}_Y$ by using the $T^1$-lifting property of $\mathrm{PD}_Y$
(cf. [Na 2], p. 281). 
But, if $\alpha \in H^2(Z, \mathbf{C})$, then it is already constructed as the objects 
of  $\mathrm{PD}^\mu_Y$ by the $T^1$-lifting property of $\mathrm{PD}^\mu_Y$. 
This means that our twistor deformation 
induces a trivial Poisson deformation of $\mu^{-1}(U)$. In particular,       
$\mathrm{Exc}(\Pi^{\alpha}_t) \ne \emptyset$ for all $t \in \mathbf{A}^1_{\alpha}$.  Q.E.D. 
\vspace{0.2cm}

Let $\mathcal{C} \subset \mathcal{X}'$ be the locus on which $\Pi$ is not an 
isomorphism.  The locus $\mathcal{C}$ can be also defined as the subset where 
$g$ is not a smooth morphism. 
We define $\mathcal{D} \subset H^2(Y, \mathbf{C})$ 
to be the closure of $g(\mathcal{C})$. By Lemma 2, $\mathcal{D}$ is a closed algebraic 
subset of $H^2(Y, \mathbf{C})$ with $\mathcal{D} \ne H^2(Y, \mathbf{C})$. 
We call $\mathcal{D}$ the {\em discriminant locus} for $\Pi$. 
\vspace{0.2cm}

{\bf Lemma 4}. {\em 
Assume that $\alpha \in H^2(Y, \mathbf{C})$ is contained in a complex hyperplane of $H^2(Y, \mathbf{C})$ generated by a codimension one face $F$ of $\overline{\mathrm{Amp}(\pi')}$ for a crepant resolution $\pi': Y' \to X$. Then $\alpha \in g(\mathcal{C})$. }  
\vspace{0.2cm}

{\em Proof}. 
The codimension one face $F$ defines a 
birational contraction map $\mu: Y' \to Z$ and it factorizes $\pi'$ as $Y' \to Z \to X$. 
Consider the universal Poisson deformation $\mathcal{Y}' \to 
H^2(Y', \mathbf{C})$ of $Y'$. We identify the base space $H^2(Y', \mathbf{C})$ 
with $H^2(Y, \mathbf{C})$. Then we have two families $\mathcal{Y}'$ and 
$\mathcal{X}'$ of Poisson varieties over 
$H^2(Y, \mathbf{C})$ and there is a birational projective morphism 
$\Pi': \mathcal{Y}' \to \mathcal{X}'$. 
By the assumption $\alpha \in \mu^*(H^2(Z, \mathbf{C}))$. By applying Lemma 3 to 
$\pi': Y' \to X$, we see that 
$\mathrm{Exc}(\Pi'_{\alpha}) \ne \emptyset$. Hence $\alpha \in g(\mathcal{C})$.  Q.E.D. 
\vspace{0.2cm}



{\bf Corollary 5}. 

{\em There are only finitely many crepant projective resolutions of $X$.}

\vspace{0.2cm}

{\em Proof}. 
Suppose that there are infinitely many crepant projective resolutions of $X$. 
Then we have infinitely many rational polyhedral cones inside $H^2(Y, \mathbf{R})$ 
corresponding to the nef cones of crepant resolutions. Note that the interior part of 
these cones do not intersect. The codimension one faces 
of these cones generates infinitely many complex hyperplanes of $H^2(Y, \mathbf{C})$. 
Since $\mathcal{D} \ne H^2(Y, \mathbf{C})$, these complex hyperplanes are all 
irreducible components of $\mathcal{D}$ by Lemma 4. This is absurd. 
Q.E.D. 
\vspace{0.2cm}

{\bf Theorem}([BCHM]). {\em Let $f: V \to W$ be a {\bf Q}-factorial terminalization of a 
normal variety $W$ with rational Gorenstein singularities. Assume that $f$ is an isomorphism in codimension one and $\rho(V/W) = 1$. Then there is a flop $f': V' \to W$, in other words, $f'$ is another {\bf Q}-factorial terminalization of $W$ such that the proper transform $D'$ of an $f$-negative divisor divisor $D$ is $f'$-ample.}    
\vspace{0.2cm}

In fact, $\bigoplus_{m \geq 0}f_*\mathcal{O}_V(mD)$ is a finitely generated $\mathcal{O}_W$-algebra. The flop $V'$ is nothing but $\mathbf{Proj}_W  \bigoplus_{m \geq 0}f_*\mathcal{O}_V(mD)$. 
\vspace{0.2cm}

{\bf Lemma 6}. {\em $\overline{\mathrm{Mov}(\pi)}$ is the union of the 
nef cones of crepant projective resolutions of $X$.}
\vspace{0.2cm}

{\em Proof}. Assume that $D$ is a $\pi$-movable divisor. 
If $D \notin \overline{\mathrm{Amp}(\pi)}$, there is a codimension one face 
of $\overline{\mathrm{Amp}(\pi)}$ such that the corresponding 
birational $X$-morphism $\phi_1: Y \to Z_1$ contracts $D$-negative curves. 
Since $D$ is $\pi$-movable, $\phi_1$ is an isomorphism in codimension one. 
By the previous theorem we can take its flop $\phi'_1 : Y_1 \to Z_1$. 
Let $\pi_1: Y_1 \to X$ be the composition of $\phi'_1$ and $Z_1 \to X$. 
Notice that $\pi_1$ is a crepant resolution by [Na 2, Corollary 31]. 
Let $D_1 \subset Y_1$ be the proper transform of $D$. If $D_1$ is $\pi_1$-nef, 
then $D \in \overline{\mathrm{Amp}(\pi_1)}$. If $D_1$ is not $\pi_1$-nef, we 
repeat the same process. Since there are only finitely many different crepant 
projective resolutions of $X$, if this process does not terminate, there appear 
the same crepant resolution twice in the process, say $Y_i$ and $Y_j$ ($i < j$). 
Take a discrete valuation $v$ of the function field $K$ of $Y_i$ in such a 
way that its center is contained in $\mathrm{Exc}(\phi_{i-1})$. Let $D_i \subset Y_i$ be 
the proper transform of $D$. Then we have inequalities for discrepancies (cf. [KMM], 
Proposition 5-1-11): 
$$ a(v, D_i) < a(v, D_{i+1}) \le a(v, D_{i+2}) \le ... \le a(v, D_j).$$ 
Here the first inequality is a strict one by the choice of $v$. 
On the other hand, we have $Y_i = Y_j$ and $D_i = D_j$ by the assumption. 
This is absurd. Thus we finally reach $\pi_n: Y_n \to X$ such that $D_n$ is 
$\pi_n$-nef. Then $D \in \overline{\mathrm{Amp}(\pi_n)}$. Q.E.D. 
\vspace{0.2cm}

The following fact is pointed out by Braden, Proudfoot and Webster.  
\vspace{0.2cm}

{\bf Proposition}([BPW, Proposition 2.19]) {\em The movable cone $\mathrm{Mov}(\pi) \subset H^2(Y, \mathbf{R})$ is a fundamental 
domain for the $W$-action on $H^2(Y, \mathbf{R})$.} 
\vspace{0.2cm}

The group $W$ also acts on $\mathcal{X}'$ and the map  $\mathcal{X}' \to H^2(Y, \mathbf{C})$ is $W$-equivarinat.
Let $\Pi': \mathcal{Y}' \to \mathcal{X}'$ be another crepant projective resolution of $\mathcal{X}'$. For $w \in W$, we take the fibre product 

\begin{equation} 
\begin{CD} 
\mathcal{Y}'_w @>>> \mathcal{Y}' \\ 
@V{\Pi'_w}VV @V{\Pi'}VV \\ 
\mathcal{X}' @>{w}>> \mathcal{X}'    
\end{CD} 
\end{equation}

We call $\mathcal{Y}'_w$ the $w$-twist of $\mathcal{Y}'$.  

Let $\{\pi_k : Y_k \to X\}_{k \in K}$ be the set of all crepant projective resolutions of 
$X$. For each $k \in K$, we can construct the universal Poisson deformation 
$\mathcal{Y}_k \to H^2(Y, \mathbf{C})$ of $Y_k$. Notice that there is a 
$H^2(Y, \mathbf{C})$-birational projective morphism $\Pi_k: \mathcal{Y}_k \to \mathcal{X}'$ and $\Pi_k$ is a crepant resolution of $\mathcal{X}'$. 
\vspace{0.2cm}

{\bf Corollary 7}. {\em  Every {\bf Q}-factorial terminalization of $\mathcal{X}'$ is 
obtained as the $w$-twist $\mathcal{Y}_{k,w}$ of $\mathcal{Y}_k$ for $w \in W$. 
In particular, they are all smooth.}
\vspace{0.2cm}

{\em Proof}. Denote by $\Pi_{k,w}$ the crepant resolution 
$\mathcal{Y}_{k,w} \to \mathcal{X}'$. 
The nef cone $\overline{\mathrm{Amp}(\Pi_{k,w})} \subset H^2(Y, \mathbf{R})$ 
is $w(\overline{\mathrm{Amp}(\pi_k)})$. By the previous proposition and Lemma 6, 
$H^2(Y, \mathbf{R})$ is covered by the cones $\{w(\overline{\mathrm{Amp}(\pi_k)})\}$. Q.E.D. 
\vspace{0.2cm}

{\bf Lemma 8}. {\em Assume that $\alpha \in \mathcal{D}$. Then there exists a 
crepant projective resolution $\Pi': \mathcal{Y}' \to \mathcal{X}'$ and a codimension 
one face of $F$ of $\overline{\mathrm{Amp}(\Pi')} \subset H^2(Y, \mathbf{R})$ 
such that the complex hyperplane $H$ of $H^2(Y, \mathbf{C})$ generated by 
$F$ contains $\alpha$.} 
\vspace{0.2cm}

{\em Proof}. We may assume that $\alpha \in g(\mathcal{C})$. Then $\Pi_{\alpha}: 
\mathcal{Y}_{\alpha} \to \mathcal{X}'_{\alpha}$ has non-empty exceptional locus\footnote{Caution: $\mathcal{Y}_{\alpha}$ is not the twistor deformation, but the 
fibre over $\alpha$}. 
Let $H$ be a $\Pi$-ample divisor. Since $\mathrm{Codim}_{\mathcal Y}
\mathrm{Exc}(\Pi) \geq 2$, we have $D: = -H \in \mathrm{Mov}(\Pi)$. Notice that 
$(D, C) < 0$ for any proper curve $C$ in $\mathcal{Y}_{\alpha}$ contracted by 
$\Pi_{\alpha}$.  Choose a codimension 
one face of $\overline{\mathrm{Amp}(\Pi)}$ and let $\phi_1 : \mathcal{Y} \to  \mathcal{Z}_1$ 
be the corresponding birational morphism. By the previous theorem there is a flop $\phi'_1: \mathcal{Y}_1 \to \mathcal{Z}_1$. We denote by $\Pi_1$ the composition of $\phi'_1$ and the 
map $\mathcal{Z}_1 \to \mathcal{X}'$. The resulting variety $\mathcal{Y}_1$ is again a 
crepant projective resolution of $\mathcal{X}'$ by Corollary 7. 
Let $D_1 \subset \mathcal{Y}_1$ be the proper transform of $D$. If $D_1$ is not $\Pi_1$-nef, 
one can find a codimension one face $F_1$ of $\overline{\mathrm{Amp}(\Pi_1)}$ in such a way that the corresponding birational map $\phi_2: \mathcal{Y}_1 \to \mathcal{Z}_2$ contracts $D_1$-negative curves. By taking the flop of $\phi_2$ we get a new crepant resolution $\Pi_2: 
\mathcal{Y}_2 \to \mathcal{X}'$. Since there are only finitely many crepant resolutions 
of $\mathcal{X}'$ by Corollary 7, this process eventually terminates by the same argument 
as in Lemma 6 and the proper transform $D_n \subset \mathcal{Y}_n$ of 
$D$ is $\Pi_n$-nef for some $n$. 
Each map $\phi_i$ induces a birational map $\phi_{i, \alpha}: \mathcal{Y}_{i-1, \alpha} \to \mathcal{Z}_{i, \alpha}$.  Suppose that all $\phi_{i, \alpha}$ are isomorphisms; then the birational map $\mathcal{Y} --\to \mathcal{Y}_n$ is an isomorphism on  a neighborhood of $\mathcal{Y}_{\alpha}$. On the other hand, $D$ is $\Pi$-negative 
around $\mathcal{Y}_{\alpha}$ and $D_n$ is $\Pi_n$-nef around $\mathcal{Y}_{n, \alpha}$. 
This is absurd. Thus  $\mathrm{Exc}(\phi_{i_0, \alpha}) \ne \emptyset$ for some $i_0$. 

By Corollary 7,  $\mathcal{Y}_{i_0-1}$ can be written as $\mathcal{Y}_{k,w}$ with a crepant resolution $Y_k$ of $X$ and $w \in W$. 
Let $C_{\alpha}$ be a proper curve on $\mathcal{Y}_{i_0-1, \alpha}$ contracted by 
$\phi_{i_0, \alpha}$. Then $([\omega_{\mathcal{Y}_{i_0-1, \alpha}}], [C_{\alpha}]) = 0$. 
By the $\mathbf{C}^*$-action on $\mathcal{Y}_{i_0-1}$ we obtain a proper curve $C$ on 
$Y_k$ as $\lim_{t \to 0} t(C_{\alpha})$. As $\mathcal{Y}_{i_0-1, \alpha}$ is  
diffeomorphic to  $Y_k$, there is a natural identification $H^2(\mathcal{Y}_{i_0-1, \alpha}, 
\mathbf{C}) \cong H^2(Y_k, \mathbf{C}) (\cong H^2(Y, \mathbf{C}))$. Under this 
identification we have $[\omega_{\mathcal{Y}_{i_0-1, \alpha}}] = \alpha$; thus   
$(\alpha, [C]) = 0$.  

Let $F$ be the codimension one face of $\overline{\mathrm{Amp}(\Pi_{i_0-1})}$ 
which determines $\phi_i$. Let $H$ be the complex hyperplane of  
$H^2(Y, \mathbf{C})$ generated by $F$. Note that $H = 
\{x \in H^2(Y, \mathbf{C}); (x, [C]) = 0\}.$ As $(\alpha, [C]) = 0$, we have 
$\alpha \in H$.  
Q.E.D. 
\vspace{0.2cm} 

{\bf Remark 9}. The locus $g(\mathcal{C})$ is invariant by $W$. By Lemma 4 and 
Corollary 7,  if $\alpha \in H^2(Y, \mathbf{C})$ is contained in the complex 
hyperplane generated by a codimension one face $F$ of $\overline{\mathrm{Amp}(\Pi')}$, 
then $\alpha \in g(\mathcal{C})$. 
Lemma 8 together with this fact implies that $\mathcal{D} = g(\mathcal{C})$ and 
that $\mathcal{D}$ coincides with the union of complex hyperplanes generated by $F
\subset \overline{\mathrm{Amp}(\Pi')}$ for all crepant projective 
resolutions $\Pi': \mathcal{Y}' \to \mathcal{X}'$ and for all codimension one faces 
$F$ of their nef cones.  
\vspace{0.2cm}

Now our results can be summarized as follows: 
\vspace{0.2cm}   
 
{\bf Main Theorem}. 

(i) {\em There are finitely many codimension one linear subspaces $\{H_i\}_{i \in I}$ of $H^2(Y, \mathbf{Q})$ such that $\mathcal{D} = \cup_i (H_i)_{\mathbf C}$. Moreover 
$\mathcal{D} = g(\mathcal{C})$. } 

(ii) {\em There are only finitely many crepant projective resolutions of $X$, say 
$\{\pi_k: Y_k \to X\}_{k \in K}$. 
The closed movable cone $\overline{\mathrm{Mov}(\pi)}$ 
is the union of nef cones: $\overline{\mathrm{Mov}(\pi)} = \bigcup_{k \in K} 
\overline{\mathrm{Amp}(\pi_k)}$.
Moreover $H^2(Y, \mathbf{R}) = \bigcup_{k \in K, w \in W}
w(\overline{\mathrm{Amp}(\pi_k)})$. }   

(iii) {\em 
The chambers $\{w(\mathrm{Amp}(\pi_k))\}_{k \in K}$ 
coincide with the open chambers determined by $\{(H_i)_{\mathbf R}\}_{i \in I}$.}  
\vspace{0.2cm}

{\bf  Remark 10}. Even if $X$ has no crepant resolution,  
Main theorem equally holds for a {\bf Q}-factorial terminalization $\pi: Y \to X$. In fact, 
we can construct a commutative diagram of universal Poisson deformations of $Y$ and $X$ (cf. [Na 3], Theorem 5.5) over affine spaces 

\begin{equation} 
\begin{CD} 
\mathcal{Y} @>>> \mathcal{X} \\ 
@V{f}VV @VVV \\ 
\mathbf{A}^n @>{q}>> \mathbf{A}^n
\end{CD} 
\end{equation}

Notice that each fibre of $\mathcal{Y} \stackrel{f}\to \mathbf{A}^n$ has singularities, 
but $f$ is a locally trivial deformation of $Y$ ([Na 3], Theorem 5.5, (b)), i.e. 
for each point $p \in Y$, there is an open subset $\mathcal{U}$ of $\mathcal{Y}$ in the classical topology such that $p \in \mathcal{U}$ 
and $f\vert_{\mathcal U}$ looks like $(\mathcal{U} \cap Y) \times \Delta^n  \to \Delta^n$.  
By the period map $p: \mathbf{A}^n \to H^2(Y, \mathbf{C})$ ([Na 5], (4.2)), the first 
affine space $\mathbf{A}^n$ is naturally identified with $H^2(Y, \mathbf{C})$. 
Moreover Theorem 1.1 of [Na 4] still holds in our case and $q$ is a finite Galois covering. 
Thus the second affine space $\mathbf{A}^n$ is identified with $H^2(Y, \mathbf{C})/W$ 
for a finite Galois group $W$. 
\vspace{0.2cm}

\begin{center}
{\bf \S 3. Examples}  
\end{center}
\vspace{0.2cm}

{\bf Example 11}.  Let $\mathfrak{g}$ be a complex simple Lie algebra and 
$G$ its adjoint group.  Let $\chi : \mathfrak{g} \to \mathfrak{g}//G$ 
be the adjoint quotient map. If $\mathfrak{h}$ is a Cartan subalgebra of 
$\mathfrak{g}$ and $W$ is the Weyl group of $\mathfrak{g}$ with respect to 
$\mathfrak{h}$, then $\mathfrak{g}//G$ can be identified with $\mathfrak{h}/W$.  
Choose a Borel subalgebra $\mathfrak{b}$ of $\mathfrak{g}$ containing $\mathfrak{h}$. 
Define a map $\mu_B: G \times^B \mathfrak{b} \to \mathfrak{g}$ by $(g, x) \to 
Ad_g(x)$. 
Let $O \subset \mathfrak{g}$ be a nilpotent orbit and $S \subset \mathfrak{g}$ 
the Slodowy slice to $O$. Consider the map $\chi\vert_S : S \to \mathfrak{h}/W$ 
and denote by $S_0$ its null fibre. By the Kostant-Kirillov form, $S_0$ is a symplectic variety 
and it admits a good $\mathbf{C}^*$-action. 
On the other hand, we define $S_B := \mu_B^{-1}(S)$ and consider the composition map 
$S_B \subset G \times^B \mathfrak{b} \to \mathfrak{h}$. We then denote by 
$S_{B,0}$ its null fibre. The map $S_B \to S$ induces a map $\pi: S_{B,0} \to S_0$, which is 
a crepant resolution of $S_0$.  We have $\mathrm{Mov}(\pi) = \mathrm{Amp}(\pi)$.  
In other words, $S_0$ has a unique crepant resolution.      
In [LNS] we proved that the commutative diagram 
 
\begin{equation} 
\begin{CD} 
S_B @>>> S \\ 
@VVV @VVV \\ 
\mathfrak{h} @>>> \mathfrak{h}/W    
\end{CD} 
\end{equation}
gives the universal Poisson deformations of $S_{B,0}$ and $S_0$ 
with some exceptions. 
The exceptional cases are when

 (i) $O$ is the subregular orbit in 
$B_n$, $C_n$ $G_2$ or $F_4$, 

(ii) $O_{[n,n]}$, $O_{[2n-2i, 2i]}$ ($1 < i \le n/2$)  in 
$C_n$  or 

(iii) 8 dimensional nilpotent orbit in $G_2$.    

When $O$ is not in the lists, we have $H^2(S_B, \mathbf{C}) \cong \mathfrak{h}$ 
and our chamber structure for $\pi$ is nothing but the Weyl chamber structure 
of $\mathfrak{h}_{\mathbf R}$.  

When $O$ is in the list (i), we know that $S_0$ is  
isomomorhic to the null fibre of the Slodowy slice to the subregular nilpotent 
orbit in another simply-laced Lie 
algebra (as symplectic $\mathbf{C}^*$-varieties).  
For $B_n$, the simply-laced Lie algebra is $A_{2n-1}$ and other cases are:  
$C_n \to D_{n+1}$, $G_2 \to D_4$, $F_4 \to E_6$. 

For (ii),  the nullfibres of the Slodowy slices for the following pairs are mutually  
isomorphic as symplectic $\mathbf{C}^*$-varieties  ([H-L], Corollary 5.2) :

($O_{[n,n]} \subset C_n$, $O_{[n+1, n+1]} \subset D_{n+1}$) 

($O_{[2n-2i, 2i]} \subset C_n$, $O_{[2n-2i+1, 2i+1]} \subset D_{n+1}$). 
\vspace{0.2cm}

Finally, the nullfibre of the Slodowy slice in the case (iii) is isomorphic 
to the null fibre of the Slodowy slice for $O_{[4, 1^2]} \subset C_3$ as symplectic  $\mathbf{C}^*$-varieties.     

In any of these cases, we have $H^2(S_B, \mathbf{C}) \cong \mathfrak{h}'$ for a Cartan subalgebra $\mathfrak{h}'$ of a suitable simple Lie algebra $\mathfrak{g}'$. Then the  chamber structure for $\pi$ again coincides with the Weyl chamber structure of $\mathfrak{h}'_{\mathbf R}$.  
\vspace{0.2cm}

{\bf Example 12}. Let $G$ be a complex simple Lie group and let $P_0$ be its 
parabolic subgroup. Let $s: T^*(G/P_0) \to \mathfrak{g}$ be the Springer map 
from the cotangent bundle of $G/P_0$ to the Lie algebra $\mathfrak{g}$ of $G$. 
The Springer map $s$ is a generically finite projective morphism. As $s$ is $G$-equivariant, 
the image of $s$ coincides with the closure $\bar{O}$ of a nilpotent orbit 
$O$ of $\mathfrak{g}$. We take the Stein factorization of $s$: $T^*(G/P_0) \stackrel{\mu}\to 
\tilde{O} \to \bar{O}$. Then $\mu$ is a crepant resolution of $\tilde{O}$. 
In [Na 6] we have studied the birational geometry related with $\mu$. According to [Na 6] 
one can describe the chamber structures of $H^2(T^*(G/P_0), \mathbf{R})$ quite explicitly
\footnote{In [Na 6] the base space $H^2(T^*(G/P_0), \mathbf{C})$ is denoted 
by $\mathfrak{k}(\mathfrak{p}_0)$.}. 
Here we shall illustrate it by using one typical example.  

We put $G = SL(5, \mathbf{C})$. Fix a maximal torus $T$ and consider the corresponding 
root system $\Phi$ of $\mathfrak{g}$. Choose a base $\{\alpha_1, \alpha_2, \alpha_3, 
\alpha_4\}$ of $\Phi$.  Denote by $\Phi^{+}$(resp. $\Phi^{-}$) the set of positive roots in $\Phi$.  
Let us consider the marked 
Dynkin diagram  

\begin{picture}(300,20)
\put(30,0){\circle{5}}\put(35,0){\line(1,0){20}}
\put(60,0){\circle*{5}}\put(65,0){\line(1,0){20}}
\put(90,0){\circle*{5}}\put(95,0){\line(1,0){20}}
\put(120,0){\circle{5}}\put(25, -9){$\alpha_1$} 
\put(55, -9){$\alpha_2$}\put(85, -9){$\alpha_3$} 
\put(115, -9){$\alpha_4$}  
\end{picture} 
\vspace{1.0cm}

Let $I$ be the set of simple roots corresponding to the white 
vertexes and denote by $\Phi_I$ the root subsystem of $\Phi$ generated by 
$I$. We set $\Phi_I^{-}:= \Phi_I \cap \Phi^{+}$. We define a parabolic subalgebra $\mathfrak{p}_0$ of $\mathfrak{g}$ by  
$$\mathfrak{p}_0 :=  \mathfrak{h} \oplus \bigoplus_{\alpha \in \Phi^{+}}\mathfrak{g}_{\alpha} \oplus \bigoplus_{\alpha \in \Phi_I^{-}} \mathfrak{g}_{\alpha}.$$  
Here $\mathfrak{h} := \mathrm{Lie}(T)$.    
Let $P_0$ be the corresponding parabolic subgroup of $G$.  Let $L \subset P_0$ be 
the Levi part of $P_0$ containing $T$. 
By Theorem 1.3 and Proposition 2.3 of [Na 6], there is a one-to-one correspondence 
between the chambers 
(inside $H^2(T^*(G/P_0), \mathbf{R})$) and the parabolic subgroups of $G$ with the Levi 
part $L$. Moreover, all such parabolic subgroups are constructed from $P_0$ by performing  
the operations ``twist"  successively.   

If we twist $P_0$ by the vertex $\alpha_2$, we have a parabolic subgroup
 $P_1$: \vspace{0.2cm}

\begin{picture}(300,20)
\put(30,0){\circle*{5}}\put(35,0){\line(1,0){20}}
\put(60,0){\circle{5}}\put(65,0){\line(1,0){20}}
\put(90,0){\circle*{5}}\put(95,0){\line(1,0){70}}
\put(170,0){\circle{5}}\put(25, -9){$-\alpha_2$} 
\put(55, -9){$-\alpha_1$}\put(85, -9){$\alpha_1 + \alpha_2 + \alpha_3$} 
\put(165, -9){$\alpha_4$}  
\end{picture} 
\vspace{1.0cm}

We can continue twisting one by one. Then we get
\vspace{0.2cm}

$P_2$ (as the twist of $P_1$ by $\alpha_1 + \alpha_2 + \alpha_3$): 

\begin{picture}(300,20)
\put(30,0){\circle*{5}}\put(35,0){\line(1,0){40}}
\put(80,0){\circle{5}}\put(85,0){\line(1,0){20}}
\put(110,0){\circle*{5}}\put(115,0){\line(1,0){90}}
\put(210,0){\circle{5}}\put(25, -9){$\alpha_3 + \alpha_4$} 
\put(75, -9){$-\alpha_4$}\put(105, -9){$-\alpha_1 - \alpha_2 - \alpha_3$} 
\put(205, -9){$\alpha_1$}  
\end{picture} 
\vspace{1.0cm}

$P_3$ (as the twist of $P_2$ by $\alpha_3 + \alpha_4$): 

\begin{picture}(300,20)
\put(30,0){\circle{5}}\put(35,0){\line(1,0){20}}
\put(60,0){\circle*{5}}\put(65,0){\line(1,0){50}}
\put(120,0){\circle*{5}}\put(125,0){\line(1,0){50}}
\put(180,0){\circle{5}}\put(25, -9){$\alpha_4$} 
\put(55, -9){$-\alpha_3-\alpha_4$}\put(115, -9){$-\alpha_1 - \alpha_2$} 
\put(175, -9){$\alpha_1$}  
\end{picture} 
\vspace{1.0cm}

$P_4$ (as the twist of $P_3$ by $-\alpha_1 - \alpha_2$): 

\begin{picture}(300,20)
\put(30,0){\circle{5}}\put(35,0){\line(1,0){20}}
\put(60,0){\circle*{5}}\put(65,0){\line(1,0){90}}
\put(160,0){\circle{5}}\put(165,0){\line(1,0){20}}
\put(190,0){\circle*{5}}\put(25, -9){$\alpha_4$} 
\put(55, -9){$-\alpha_2 -\alpha_3 -\alpha_4$}\put(155, -9){$-\alpha_1$} 
\put(185, -9){$\alpha_1+\alpha_2$}  
\end{picture} 
\vspace{1.0cm}

$P_5$ (as the twist of $P_4$ by $\alpha_2 + \alpha_3 + \alpha_4$): 

\begin{picture}(300,20)
\put(30,0){\circle{5}}\put(35,0){\line(1,0){20}}
\put(60,0){\circle*{5}}\put(65,0){\line(1,0){90}}
\put(160,0){\circle{5}}\put(165,0){\line(1,0){30}}
\put(200,0){\circle*{5}}\put(25, -9){$\alpha_1$} 
\put(55, -9){$\alpha_2 + \alpha_3 + \alpha_4$}
\put(155, -9){$-\alpha_4$} 
\put(195, -9){$-\alpha_3$}  
\end{picture} 
\vspace{1.0cm}

These 6 parabolic subgroups are the complete lists of those with the Levi part $L$. 
In the remainder we put $Y := T^*(G/P_0)$, $\mathcal{Y} := G \times^{P_0}r(\mathfrak{p}_0)$ 
and $\mathcal{Y}_i := G \times^{P_i}r(\mathfrak{p}_i)$. Here $r(\mathfrak{p}_i)$ are solvable 
radicals of $\mathfrak{p}_i$.   
Since $Y$ is a vector bundle over $G/P$, we have an isomorphism 
$H^2(Y, \mathbf{R}) \cong H^2(G/P, \mathbf{R})$. On the other hand, we have 
$$H^2(G/P, \mathbf{R}) \cong \mathrm{Hom}_{alg.gp}(L, \mathbf{C}^*)\otimes_{\mathbf Z}
\mathbf{R} = \{x \in \mathfrak{h}^*_{\mathbf R}; (x, \alpha_1) = (x, \alpha_4) = 0\},$$
where $(\;, \;)$ is the Killing form on $\mathfrak{h}^*_{\mathbf R}$. 
Thus $H^2(Y, \mathbf{R})$ can be identified with the dual space of the 2-dimensional 
vector space $\mathbf{R}\alpha_2 \oplus \mathbf{R}\alpha_3$ by the Killing form.  
Now $H^2(Y, \mathbf{R})$ is decomposed into 6 chambers:

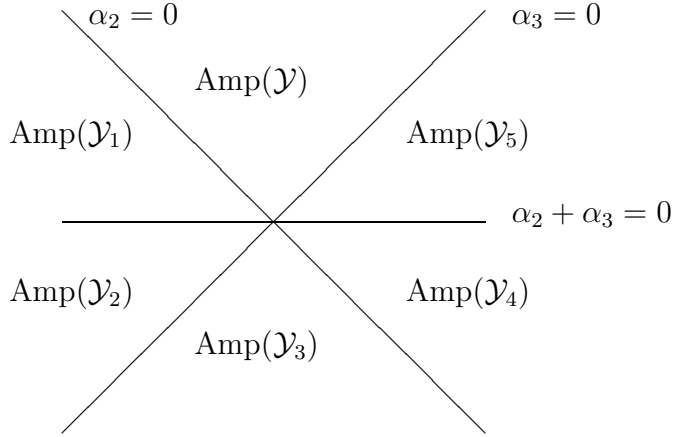
\begin{picture}(300,100)(0,0) 
\put(150,0){\line(1,0){80}}\put(150,0){\line(1,1){80}}\put(150,0){\line(-1,1){80}}
\put(150,0){\line(-1,0){80}}
\put(150,0){\line(-1,-1){80}}
\put(150,0){\line(1,-1){80}}
\put(120, 50){$\mathrm{Amp}(\mathcal{Y})$}\put(200, 30){$\mathrm{Amp}(\mathcal{Y}_5)$}
\put(200,-30){$\mathrm{Amp}(\mathcal{Y}_4)$}\put(120,-50){$\mathrm{Amp}(\mathcal{Y}_3)$} 
\put(50,-30){$\mathrm{Amp}(\mathcal{Y}_2)$}\put(50,30){$\mathrm{Amp}(\mathcal{Y}_1)$}
\put(80,75){$\alpha_2 = 0$}\put(240,75){$\alpha_3 = 0$}\put(240,0){$\alpha_2 + \alpha_3 
= 0$}
\end{picture} 
\vspace{5.0cm}
 
The movable cone $\mathrm{Mov}(Y)$ is the upper half plane. The Weyl group $W$ 
is isomorphic to $\mathbf{Z}/2\mathbf{Z}$, which acts on $H^2(Y, \mathbf{R})$ as 
the reflection with respect to the line $\{\alpha_2 + \alpha_3 = 0\}$.   
\vspace{0.2cm}

{\bf Example 13}. The following example is due to Ryo Yamagishi.  It clarifies a  
difference between symplectic singularities and non-symplectic singularities. 

Let us consider the line bundle $Z := p_1^*\mathcal{O}_{{\mathbf P}^1}(-1) 
\otimes p_2^*\mathcal{O}_{{\mathbf P}^1}(-1) \otimes p_3^*\mathcal{O}_{{\mathbf P}^1}(-1)$  on $\mathbf{P}^1 \times \mathbf{P}^1 
\times \mathbf{P}^1$. Denote by $\Sigma \subset Z$ the zero section of the bundle map 
$Z \to \mathbf{P}^1 \times \mathbf{P}^1 \times \mathbf{P}^1$. For $1 \le i < j \le 3$, let 
$p_{ij}: \Sigma \to \mathbf{P}^1 \times \mathbf{P}^1$ be the projection to the $i$-th and $j$-th factors. There is a birational projective morphism $f_{ij}: Z \to X_{ij}$ to a smooth variety 
$X_{ij}$ such that $f_{ij}\vert_{\Sigma} = p_{ij}$ and $f_{ij}$ is an isomorphism outside $\Sigma$.  
We put $\Sigma_{ij} := p_{ij}(\Sigma)$. Then $\Sigma_{ij} \cong \mathbf{P}^1 \times \mathbf{P}^1$ and $f_{ij}$ is the blowing up of $X_{ij}$ along $\Sigma_{ij}$. Moreover there 
is a birational projective morphism $f: Z \to X$ to an affine (singular) variety $X$, which 
contracts $\Sigma$ to a point. By the construction $f$ factors through $f_{ij}$ for each pair 
$(i,j)$ with $1 \le i < j \le 3$.  Then $X_{12}$, $X_{13}$ and $X_{23}$ are three different crepant 
resolutions of $X$. 

Recall that $\Sigma_{ij}$ is naturally identified with the product $\mathbf{P}^1 \times \mathbf{P}^1$ of the $i$-th and the $j$-th factors 
of $\mathbf{P}^1 \times \mathbf{P}^1 \times \mathbf{P}^1$. Let $l_{ij}^i$ be a fibre of the projection map $q_{ij}^i: \Sigma_{ij} \to \mathbf{P}^1$ to the $i$-th factor.. Then the normal bundle of such a curve (inside $X_{ij}$) is isomorphic to 
$\mathcal{O}_{\mathbf{P}^1}(-1)^{\oplus 2} \oplus \mathcal{O}_{\mathbf{P}^1}$.
 
Now one has a birational projective morphism $g_{ij}^i: X_{ij} \to X_{ij}^i$ such that 
$g_{ij}^i\vert_{\Sigma_{ij}} = q_{ij}^i$ and $g_{ij}^i$ is an isomorphism outside $\Sigma_{ij}$. 
By the definition $X_{12}^1 = X_{13}^1$, $X_{12}^2 = X_{23}^2$ and $X_{13}^3 = X_{23}^3$. 
Thus we shall denote them by $X_1$, $X_2$ and $X_3$. Note that these 4-folds $X_i$ have  
3-dimensional ordinary double points along $\mathrm{Sing}(X_i) (\cong \mathbf{P}^1)$. The 4-folds $X_{12}$, $X_{13}$ 
and $X_{23}$ are related by flops
 
$$ X_{12} \to X_2 \leftarrow X_{23} $$
$$ X_{23} \to X_3 \leftarrow X_{13}$$ 
$$ X_{13} \to X_1 \leftarrow X_{12}.$$               

Then $\mathrm{Pic}(X_{12}) \otimes_{\mathbf{Z}}\mathbf{R}$ is divided by 
the ample cones of $X_{12}$, $X_{23}$ and $X_{13}$. There are 3 open chambers and they 
are determined by half-lines. Thus $X$ is not a symplectic singularity.

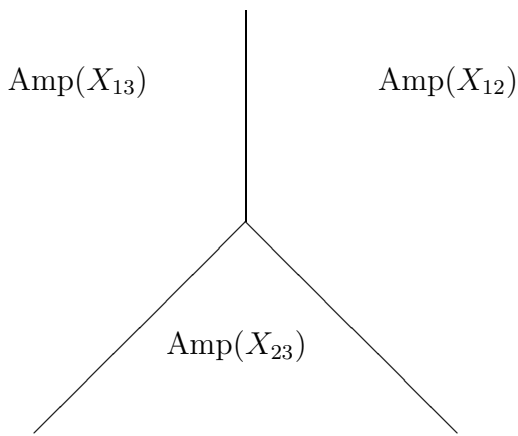
\begin{picture}(300,100)(0,0) 
\put(150,0){\line(0,1){80}}\put(150,0){\line(1,-1){80}}
\put(150,0){\line(-1,-1){80}}
\put(200, 50){$\mathrm{Amp}(X_{12})$}\put(60, 50){$\mathrm{Amp}(X_{13})$}
\put(120,-50){$\mathrm{Amp}(X_{23})$}
\end{picture} 
\vspace{5.0cm}

\begin{center}
{\bf References} 
\end{center}

[A-W] Andreatta, M., Wisniewski, J.: 
4-dimensional symplectic contractions, to appear in Geom. Dedicata 

[BCHM] Birkar, C., Cascini, P., Hacon, C., McKernan, J.: 
Existence of minimal models for varieties of log general type, J. Amer. Math. Soc. {\bf 23} (2010), 405-468 

[Be] Beauville, A.: Symplectic singularities, Invent. Math. {\bf 139} (2000),  
541-549

[BPW] Braden, T., Proudfoot, N., Webster, B.: 
Quantizations of conical symplectic resolutions I: local and global structures, arXiv: 1208.3863 

[G-K] Ginzburg, V., Kaledin, D.: 
Poisson deformations of symplectic quotient singularities, Adv. Math. {\bf 186} (2004), 
1-57 

[H-K] Hu, Y., Keel, S.: 
Mori dream spaces and GIT, Michigan Math. J. {\bf 48} (2000), 331-348 

[H-L] Henderson, A., Licata, A.: Diagonal automorphisms of 
quiver varieties, arXiv: 1309.0572

[Ka] Kaledin, D.: Symplectic resolutions: deformations and birational maps, 
math/0012008 

[K-V] Kaledin, D., Verbitsky, M.: Period maps for non-compact holomorphically 
symplectic manifolds, Geom. Funct. Anal. {\bf 12} (2002), 1265-1295

[KMM] Kawamata, Y., Matsuda, K., Matsuki, K.: 
Introduction to the minimal model problem , Adv. Stud. Pure Math. {\bf 10} (1985), 
283-360  

[LNS] Lehn, M., Namikawa, Y., Sorger, C.: 
Slodowy slices and universal Poisson deformations, Compos. Math. {\bf 148} 
(2012), 121-144  
 
[Na 1] Namikawa, Y.: Equivalence of symplectic singularities, to appear in Kyoto J. Math. 
{\bf 53}, no. 2, 2013 (Memorial issue of Prof. Maruyama) 

[Na 2] Namikawa, Y.: Flops and Poisson deformations of symplectic varieties, 
Publ. Res. Inst. Math. Sci. {\bf 44} (2008), 259-314 

[Na 3] Namikawa, Y.: Poisson deformations of affine varieties, 
Duke Math. J. {\bf 156} (2011), 51-85 

[Na 4] Namikawa, Y.: Poisson deformations of affine symplectic varieties, II, 
Kyoto J. Math. {\bf 50} (2010), 727-752 
 
[Na 5] Namikawa, Y.: Birational geometry for nilpotent orbits, in Handbook of 
Moduli (ed. G. Farkas, I. Morrison), vol. III, Advanced Lectures in Mathematics {\bf 26} 
(2013), 1-38  

[Na 6] Namikawa, Y.: Birational geometry and deformations of nilpotent 
orbits, Duke Math. J. {\bf 143} (2008), 375-405

[Sch] Schlessinger, M.: Functors of Artin rings, Trans. Amer. Math. Soc. {\bf 130} (1968), 
208-222

[Wa] Wahl, J.: Equisingular deformations of normal surface singularities, 
I, Ann. of Math. {\bf 104} (1976), 325-356

\vspace{0.2cm}

\begin{center}
Department of Mathematics, Faculty of Science, Kyoto University, Japan 

namikawa@math.kyoto-u.ac.jp
\end{center}

\end{document}